\documentclass{article}
\usepackage{amssymb}
\usepackage{amsmath}

\newtheorem{theorem}{Theorem}

\newtheorem{definition}[theorem]{Definition}

\newtheorem{lemma}[theorem]{Lemma}

\newtheorem{proposition}[theorem]{Proposition}
\newtheorem{remark}[theorem]{Remark}

\input{tcilatex}
\begin{document}

\title{Hylomorphic solitons on lattices.}
\author{Vieri Benci$^{\ast}$, Donato Fortunato$^{\ast\ast}$ \\
%EndAName
$^{\ast}$Dipartimento di Matematica Applicata \textquotedblleft U.
Dini\textquotedblright\\
Universit\`{a} di Pisa\\
via F. Buonarroti 1/c 56127 Pisa, Italy\\
e-mail: benci@dma.unipi.it\\
$^{\ast\ast}$Dipartimento di Matematica \\
Universit\`{a} di Bari and INFN sezione di Bari\\
Via Orabona 4, 70125 Bari, Italy\\
e-mail: fortunat@dm.uniba.it}
\maketitle

\textit{To our friend and colleague Louis Nirenberg, with affection and
admiration.}

\tableofcontents

\bigskip

\


\bigskip

\bigskip

\bigskip

\section{Introduction}

Roughly speaking a \textit{solitary wave} is a solution of a field equation
whose energy travels as a localized packet and which preserves this
localization in time. A \textit{soliton} is a solitary wave which exhibits
some strong form of stability so that it has a particle-like behavior (see
e.g. \cite{Ba-Be-R.}, \cite{milano}, \cite{befogranas}, \cite{raj}, \cite%
{vil}, \cite{yangL}).\

To day, we know (at least) three mechanisms which might produce solitary
waves, vortices and solitons:

\begin{itemize}
\item Complete integrability, (e.g. Kortewg-de Vries equation);

\item Topological constraints, (e.g. Sine-Gordon equation);

\item Ratio energy/charge: (e.g. the nonlinear Schroedinger equation (NS)
and the nonlinear Klein-Gordon equation (NKG)).
\end{itemize}

Following \cite{hylo}, \cite{BBBM}, \cite{milano}, the third type
of solitary waves or solitons will be called \textit{hylomorphic}.
This class includes the $Q$-\textit{balls }which are spherically
symmetric solutions of NKG and vortices which might be defined as
\textit{spinning} Q-balls. Also it includes solitary waves and
vortices which occur, by the same mechanism, in NS and in gauge
theories; a bibliography on this subject can be found in the
review papers \cite{befogranas}, \cite{milano}, and, for the
vortices, in \cite{befo10}.

This paper is devoted to the proof of a general abstract theorem which can
be applied to the main situations considered in the literature (see e.g.
\cite{milano} and \cite{befogranas}). However this theorem can be also
applied to the NS and to the NKG defined on a lattice (namely eq. (\ref{NSV}%
) when $V$ satisfies (\ref{vii}) and eq. (\ref{KG}) when $W$ satisfies (\ref%
{gii})). These results are new.

The paper is organized as follows. In section \ref{HS} we give the
definition of hylomorphic solitons and describe their general features; in
section \ref{AR} we prove some abstract results on the existence of
hylomorphic solitons; in section \ref{SSE} and in section \ref{SKG} we apply
the abstract results to NS and to NKG defined on a lattice.

\section{Hylomorphic solitary waves and solitons\label{HS}}

\subsection{An abstract definition of solitary waves and solitons\label{be}}

Solitary waves and solitons are particular \textit{orbits} of a dynamical
system described by one or more partial differential equations. The states
of this system are described by one or more \textit{fields} which
mathematically are represented by functions
\begin{equation}
u:\mathbb{R}^{N}\rightarrow V  \label{lilla}
\end{equation}%
where $V$ is a vector space with norm $\left\vert \ \cdot \ \right\vert _{V}$
which is called the internal parameters space. We denote our dynamical
system by $\left( X,T\right) $ where $X$ is the set of the states and $T:%
\mathbb{R}\times X\rightarrow X$ is the time evolution map. If $u_{0}\in X,$
the evolution of the system will be described by the function
\begin{equation}
U\left( t,x\right) =T_{t}u_{0}(x),\text{ }t\in \mathbb{R}\text{, }x\in
\mathbb{R}^{N}.  \label{flusso}
\end{equation}

Now we can give a formal definition of solitary wave:

\begin{definition}
\label{solw}An orbit $U\left( t,x\right) $ is called solitary wave if it has
the following form:%
\begin{equation*}
U\left( t,x\right) =h(t,x)u_{0}(\gamma \left( t\right) x)
\end{equation*}%
where
\begin{equation*}
\gamma \left( t\right) :\mathbb{R}^{N}\rightarrow \mathbb{R}^{N}
\end{equation*}%
is a one parameter group of isometries which depends smoothly on $t$ and%
\begin{equation*}
h(t,x):V\rightarrow V
\end{equation*}%
is a group of (linear) transformation of the internal parameter space which
depends smoothly on $t$ and $x$. In particular, if $\gamma \left( t\right)
x=x,$ $U$ is called standing wave.
\end{definition}

For example, consider a solution of a field equation which has the following
form
\begin{equation}
U\left( t,x\right) =u_{0}(x-x_{0}-\mathbf{v}t)e^{i(\mathbf{v\cdot x-}\omega
t)};\ u_{0}\in L^{2}(\mathbb{R}^{N},\mathbb{C});  \label{mer}
\end{equation}%
In this case%
\begin{eqnarray*}
\gamma \left( t\right) x &=&x-x_{0}-\mathbf{v}t \\
h(t,x) &=&e^{i(\mathbf{v\cdot x-}\omega t)}.
\end{eqnarray*}%
In this paper we are interested in standing waves, so (\ref{mer}) takes the
form%
\begin{equation}
U\left( t,x\right) =u_{0}(x-x_{0})e^{-i\omega t};\ u_{0}\in L^{2}(\mathbb{R}%
^{N},\mathbb{C});  \label{mers}
\end{equation}%
The solitons are solitary waves characterized by some form of stability. To
define them at this level of abstractness, we need to recall some well known
notions from the theory of dynamical systems.

A set $\Gamma \subset X$ is called \textit{invariant} if $u_{0}\in \Gamma
\Rightarrow \forall t,T_{t}u_{0}\in $ $\Gamma ;$ an invariant set $\Gamma
\subset X$ is called \textit{isolated} if it has a neighborhood $N$ such
that:
\begin{equation*}
\text{if }\Delta \ \text{is\ an\ invariant\ set\ and\ }\Gamma \subset \Delta
\subset N,\text{ then }\Gamma =\Delta .
\end{equation*}

\begin{definition}
Let $\left( X,d\right) $ be a metric space and let $\left( X,T\right) $ be a
dynamical system. An isolated invariant set $\Gamma \subset X$ is called
stable, if $\forall \varepsilon >0,$ $\exists \delta >0,\;\forall u\in X$,
\begin{equation*}
d(u,\Gamma )\leq \delta ,
\end{equation*}%
implies that
\begin{equation*}
\forall t\in \mathbb{R},\;\;d(T_{t}u,\Gamma )\leq \varepsilon .
\end{equation*}
\end{definition}

Now we are ready to give the definition of (standing) soliton:\bigskip

\begin{definition}
\label{ds} A standing wave $U(t,x)$ is called soliton if there is an
isolated invariant set $\Gamma $ such that

\begin{itemize}
\item (i) $\forall t,\ U(t,\cdot )\in \Gamma $

\item (ii) $\Gamma $ is stable

\item (iii) $\Gamma $ has the following structure%
\begin{equation}
\Gamma =\left\{ u(x-x_{0}):u\in \mathcal{K},\ x_{0}\in F\right\} \cong
\mathcal{K}\times F  \label{is}
\end{equation}%
where $\mathcal{K}\subset X$ is a compact set and $F\subseteq \mathbb{R}^{N}$
is not necessarily compact.
\end{itemize}
\end{definition}

\bigskip

A generic field $u\in \Gamma $ can be written as follows%
\begin{equation*}
u_{\theta }(x-x_{0}),
\end{equation*}%
where $\theta $ belongs to a set of indices $\Xi $ which parametrize $%
\mathcal{K}$\ and $x_{0}\in F.$

In the generic case of many concrete situations, $\Gamma $ is a manifold,
then, (iii) implies that it is finite dimensional and $\left( \theta
,x_{0}\right) $ are a system of coordinates.

For example, in the case (\ref{mer}), we have%
\begin{eqnarray*}
\mathcal{K} &=&\left\{ u_{0}(x)e^{i\theta }:\theta \in \mathbb{R}/\left(
2\pi \mathbb{Z}\right) \right\} \\
\Gamma &=&\left\{ u_{0}(x-x_{0})e^{i\theta }:\theta \in \mathbb{R}/\left(
2\pi \mathbb{Z}\right) ,\ x_{0}\in \mathbb{R}^{N}\right\} \cong \mathcal{K}%
\times \mathbb{R}^{N} \\
\gamma \left( t\right) x &=&x \\
\theta \left( t\right) &=&\mathbf{-}\omega t;
\end{eqnarray*}%
then $U_{0}(t,x)=u_{0}(x-x_{0})e^{-i\omega t}\in \Gamma $ is a soliton if $%
\Gamma $ is stable.

The proof of (ii) of definition \ref{ds}, namely that a solitary wave has
enough stability to be a soliton, is a delicate question and in many cases
of interest it is open. Moreover the notion of stability depends on the
choice of the space $X$ and on the choice of its metric $d$ and hence
different choices might lead to more or less satisfactory answers.

\bigskip

\subsection{Integrals of motion and hylomorphic solitons\label{im}}

The existence and the properties of hylomorphic solitons are guaranteed by
the interplay between \textit{energy }$E$ and another integral of motion
which, in the general case, is called \textit{hylenic charge} and it will be
denoted by\textit{\ }$C.$ Notice that $E$ and $C$ can be considered as
functionals defined on the phase space $X.$

Thus, we make the following assumptions on the dynamical system $(X,T):$

\begin{itemize}
\item \textbf{A-1. }\textit{It is variational, namely the equations of the
motions are the Euler-Lagrange equations relative to a Lagrangian density $%
\mathcal{L}$}$\left[ u\right] $\textit{.}

\item \textbf{A-2. }\textit{The equations are invariant for time
translations. }

\item \textbf{A-3. }\textit{The equations are invariant for a }$S^{1}$-%
\textit{action acting on the internal parameter space }$V$ (cf. (\ref{lilla}%
)).
\end{itemize}

By \textbf{A-1,\ A-2} and \textbf{A-3} and Noether's theorem (see e.g. \cite%
{Gelfand}, \cite{milano}, \cite{befogranas}) it follows that our dynamical
system has 2 first integrals:

\begin{itemize}
\item the invariance with respect to time translations gives the
conservation of the energy which we will shall denote by $E(u)$;

\item The invariance \textbf{A-3} gives another integral of motion called
\textit{hylenic charge} which we shall denote by $C(u)$.
\end{itemize}

Now we set%
\begin{equation*}
\mathfrak{M}_{\sigma }=\left\{ u\in X:\ \left\vert C(u)\right\vert =\sigma
\right\} .
\end{equation*}

Using the definition \ref{ds}, we get the definition of hylomorphic soliton
as follows:

\begin{definition}
\label{dss}Let $U$ be a soliton according to definition \ref{ds}. $U$ is
called a (standing) \emph{hylomorphic soliton } if $\Gamma $ (as defined in (%
\ref{is})) coincides with the set of minima of $E\ $on $\mathfrak{M}_{\sigma
},$ namely
\begin{equation*}
\Gamma _{\sigma }=\left\{ u\in \mathfrak{M}_{\sigma }:E(u)=c_{\sigma
}\right\}
\end{equation*}%
with%
\begin{equation*}
c_{\sigma }=\ \underset{u\in \mathfrak{M}_{\sigma }}{\min }E(u).
\end{equation*}
\end{definition}

\bigskip

\begin{remark}
Suppose that the Lagrangian \textit{$\mathcal{L}$}$\left[ u\right] $ is
invariant for (a representation of) the Lorentz or the Galileo group. Then
given a standing hylomorphic soliton, we can get a \emph{hylomorphic
travelling soliton }just by Galileo or a Lorentz transformation
respectively.\bigskip
\end{remark}

\bigskip

We recall that in physics literature the solitons of definition \ref{dss}
are called Q-balls \cite{Coleman86} and were first studied in the pioneering
paper \cite{rosen68}. The existence of stable solitary waves in particular
cases has been extablished in \cite{CL82} and \cite{shatah}. The existence
of hylomorphic solitons in more general  equations has been proved in \cite%
{BBBM}.

If the energy $E$ is unbounded from below on $\mathfrak{M}_{\sigma }$ it is
still possible (\cite{strauss}, \cite{Beres-Lions}) to have standing wave
(see def. \ref{solw}). Moreover there are also cases \cite{shatah} in which
it is possible to have solitons (see def. \ref{ds}) which are only local
minimizers \cite{bonanno}. These solitons are not hylomorphic (def. \ref{dss}%
) and they can be destroyed by a perturbation which send them out of the
basin of attraction. \textsc{\ }

In the next section we analyze some abstract situations which imply $\Gamma
_{\sigma }\neq \varnothing $ and the existence of \emph{hylomorphic solitons}
(definition \ref{dss}).

\section{Abstract results\label{AR}}

\subsection{The general framework\label{gf}}

We assume that $E$\ and $C$\ are two functionals on $\mathcal{D}\left(
\mathbb{R}^{N},V\right) (\equiv C_{0}^{\infty }(\mathbb{R}^{N}))$ defined by
densities.\textit{\ }This means that, given $u\in \mathcal{D}\left( \mathbb{R%
}^{N},V\right) ,$ there exist density functions $\rho _{E,u}\left( x\right) $
and $\rho _{C,u}\left( x\right) \in L^{1}(\mathbb{R}^{N})$ i. e. functions
such that%
\begin{eqnarray*}
E\left( u\right) &=&\int \rho _{E,u}\left( x\right) \ dx \\
C\left( u\right) &=&\int \rho _{C,u}\left( x\right) \ dx.
\end{eqnarray*}%
Also we assume that the energy can be written as follows
\begin{equation*}
E\left( u\right) =\frac{1}{2}\int \rho _{E,u}^{(2)}\left( x\right) \ dx+\int
\rho _{E,u}^{(3)}\left( x\right) \ dx
\end{equation*}%
where $\rho _{E,u}^{(2)}$ is quadratic in $u$ and $\rho _{E,u}^{(3)}$
contains the higher order terms.

If we assume $\rho _{E,u}^{(2)}>0$ for $u\neq 0,$ then we can define the
following norm:%
\begin{equation}
\left\Vert u\right\Vert ^{2}=\int \rho _{E,u}^{(2)}\left( x\right) \ dx
\label{norm}
\end{equation}%
and the Hilbert space%
\begin{equation*}
X=\left\{ \text{closure of }\mathcal{D}\left( \mathbb{R}^{N},V\right) \
\text{with respect to }\left\Vert u\right\Vert \right\} .
\end{equation*}

We assume that the energy $E$ and the charge $C$ can be extended as
functional of class $C^{2}$ in $X;$ in particular we will write $E$ as
follows:

\begin{equation}
E\left( u\right) =\frac{1}{2}\left\langle Lu,u\right\rangle +K(u)
\label{elle}
\end{equation}%
where $L:X\rightarrow X^{\prime }$ is the duality operator, namely $%
\left\langle Lu,u\right\rangle =\left\Vert u\right\Vert ^{2}$ and $K$ is
superquadratic. Also, we assume that%
\begin{equation*}
C\left( 0\right) =0;\ C^{\prime }\left( 0\right) =0.
\end{equation*}%
so that we can write
\begin{equation}
C(u)=\left\langle L_{0}u,u\right\rangle +K_{0}(u)  \label{cu}
\end{equation}%
where $L_{0}$ is a linear operator and $K_{0}$ is superquadratic.

For any $\Omega \subset \mathbb{R}^{N}$ we will write
\begin{eqnarray*}
E_{\Omega }\left( u\right) &=&\int_{\Omega }\rho _{E,u}\left( x\right) \ dx
\\
C_{\Omega }\left( u\right) &=&\int_{\Omega }\rho _{C,u}\left( x\right) \ dx
\\
\left\Vert u\right\Vert _{\Omega }^{2} &=&\int_{\Omega }\rho
_{E,u}^{(2)}\left( x\right) \ dx \\
K_{\Omega }\left( u\right) &=&\int_{\Omega }\rho _{E,u}^{(3)}\left( x\right)
\ dx.
\end{eqnarray*}

In the general scheme described above, we make the following assumptions:

\begin{itemize}
\item (E-0) (\textbf{the main protagonists}) $E$\ and $C$\ are two
functionals on $X$ of the form (\ref{elle}) and (\ref{cu})

\item (E-1) (\textbf{lattice translation invariance}) \textit{the charge and
the energy are lattice translation invariant.}
\end{itemize}

Namely we have that $\forall z\in \mathbb{Z}^{N}$%
\begin{eqnarray*}
E\left( T_{z}u\right) &=&E\left( u\right) \\
C\left( T_{z}u\right) &=&C\left( u\right)
\end{eqnarray*}%
where $T_{z}:X\rightarrow X\ $is a linear representation of the additive
group $\mathbb{Z}^{N}$ defined as follows:

\begin{equation}
u(x)=u(x+Az)  \label{lt}
\end{equation}%
$A$ is an invertible matrix which characterizes the representation $T_{z}.$
Such a $T_{z}$ will be called \textit{lattice transformation}\textbf{.}

\begin{itemize}
\item (E-2) \textbf{(coercivity})\textit{\ if} $E\left( u_{n}\right) $
\textit{and} $C\left( u_{n}\right) $ \textit{are bounded, then} $\left\Vert
u_{n}\right\Vert $ \textit{is bounded}

\item (E-3) \textbf{(local compactness) }\textit{namely, if} $%
u_{n}\rightharpoonup \bar{u},$ \textit{weakly in} $X,$ \textit{then for
bounded} $\Omega $
\begin{eqnarray*}
K_{\Omega }\left( u_{n}\right) &\rightarrow &K_{\Omega }\left( \bar{u}\right)
\\
C_{\Omega }\left( u_{n}\right) &\rightarrow &C_{\Omega }\left( \bar{u}\right)
\end{eqnarray*}

\item (E-4) \textbf{(boundedness) }\textit{if} $\left\Vert u\right\Vert \leq
M,$ \textit{then} $K_{\Omega }^{\prime }\left( u\right) \ $and$\ C_{\Omega
}^{\prime }(u)$ \textit{are bounded in} $X^{^{\prime }}\left( \Omega \right)
$ \textit{for }\emph{any} $\Omega \subset \mathbb{R}^{N}$
\end{itemize}

\subsection{The main theorems\label{sat}}

In the framework of the previous section, we want to investigate sufficient
conditions which guarantee that the energy has a minimum on the set

\begin{equation*}
\mathfrak{M}_{\sigma }=\left\{ u\in X:\left\vert C\left( u\right)
\right\vert =\sigma \right\} ,
\end{equation*}%
namely that $\Gamma _{\sigma }\neq \varnothing $ where%
\begin{equation*}
\Gamma _{\sigma }=\left\{ u\in \mathfrak{M}_{\sigma }:E(u)=\ \underset{v\in
\mathfrak{M}_{\sigma }}{\min }E(v)\right\} .
\end{equation*}%
In this section and in the next one we will study this minimization problem,
namely we may think of $E$ and $C$ as two abstract functionals. In section %
\ref{cmt} we will apply the minimization result to the case in which $E$ and
$C$ are just the energy and the hylenic charge of a dynamical system.

We set%
\begin{eqnarray}
Q_{0} &=&\left\{ x=(x_{1},...,x_{N})\in \mathbb{R}^{N}:0\leq x_{i}<1,\
i=1,..,N\right\} \ \ \ and  \label{qu0} \\
Q &=&AQ_{0}  \label{qu}
\end{eqnarray}%
where $A$ is the matrix in (\ref{lt}). Also we set%
\begin{equation}
e_{0}=\ \underset{\varepsilon \rightarrow 0}{\lim }\ \inf \ \left\{ \frac{%
E_{Q}(u)}{\left\vert C_{Q}(u)\right\vert }:\text{\ }u\in X,\ \left\Vert
u\right\Vert _{Q}\leq \varepsilon ,\ \left\vert C_{Q}(u)\right\vert
>0\right\} .  \label{new}
\end{equation}

The value $e_{0}=+\infty $ is allowed.

We now set,

\begin{equation}
\Lambda \left( u\right) =\ \frac{E\left( u\right) }{\left\vert C\left(
u\right) \right\vert }  \label{landa}
\end{equation}

\begin{equation}
\Lambda _{\sigma }=\ \underset{\varepsilon \rightarrow 0}{\lim }\ \inf
\left\{ \Lambda \left( u\right) :\sigma -\varepsilon \leq \left\vert C\left(
u\right) \right\vert \leq \sigma +\varepsilon \right\} .  \label{nuovo}
\end{equation}

\begin{equation}
\Lambda ^{\ast }=\ \underset{u\in X}{\inf }\Lambda \left( u\right) .
\label{brutta1}
\end{equation}

\bigskip

\begin{theorem}
\label{due} Assume (E-0,..,E-4). Moreover assume that%
\begin{equation}
0<\Lambda ^{\ast }<e_{0}.  \label{no}
\end{equation}%
Then, there exists $\bar{\sigma}$ such that%
\begin{equation}
\Gamma _{_{\bar{\sigma}}}\neq \varnothing  \label{sfigata}
\end{equation}%
where $\Gamma _{_{\bar{\sigma}}}$ is as in definition \ref{dss}. Moreover

\begin{itemize}
\item if $u_{n}$ is a sequence such that $\Lambda \left( u_{n}\right)
\rightarrow \Lambda _{_{\bar{\sigma}}}$ and $\left\vert C\left( u_{n}\right)
\right\vert \rightarrow \bar{\sigma}$ then
\begin{equation}
dist(u_{n},\Gamma _{_{\bar{\sigma}}})\rightarrow 0  \label{sfigata0}
\end{equation}

\item $\Gamma _{_{\bar{\sigma}}}$ has the structure in (\ref{is}), namely%
\begin{equation}
\Gamma _{_{\bar{\sigma}}}=\left\{ u(x-x_{0}):u(x)\in \mathcal{K},\ x_{0}\in
F\right\}  \label{sfigata1}
\end{equation}%
with $\mathcal{K}$ compact and $F\subset \mathbb{R}^{N}$ is a closed set
such that $F=j+F,\ \forall j\in \mathbb{Z}^{N}.$

\item any $u\in \Gamma _{_{\bar{\sigma}}}$ solves the eigenvalue problem%
\begin{equation}
E^{\prime }(u)=\lambda C^{\prime }(u).  \label{sfigata2}
\end{equation}
\end{itemize}
\end{theorem}

In many concrete situations $C(u)\ $and $\Lambda \left( u\right) $ behave
monotonically with respect to the action of the dilatation group $R_{\theta
} $ defined by%
\begin{equation*}
R_{\theta }u(x)=u(\theta x),\ \theta \in \mathbb{R}^{+}.
\end{equation*}%
In this case we obtain a stronger result:

\begin{theorem}
\label{tre} Let the assumptions of Th. \ref{due} be satisfied. Moreover
suppose that there is the action of a group $R_{\theta },\ \theta \in
\mathbb{R}^{+}\ $such that $\Lambda \left( R_{\theta }u\right) $ is
decreasing in $\theta \ $while $E(R_{\theta }u)$ and $C(R_{\theta }u)\ $are\
increasing. Then there exists $\sigma _{0}$ such that, for any $\bar{\sigma}%
\geq \sigma _{0},$ the same conclusions of Th. \ref{due} hold$.$
\end{theorem}

The following proposition gives an expression of $e_{0}$ (see (\ref{new}))
which will be useful in the applications of Theorems \ref{due} and \ref{tre}.

\begin{proposition}
\label{laida} Let $E$ and $C$ be as in (\ref{elle}) and (\ref{cu}); then
\begin{equation}
e_{0}=\ \underset{u\in X}{\inf }\ \frac{\frac{1}{2}\left\langle
Lu,u\right\rangle _{Q}}{\left\vert \left\langle L_{0}u,u\right\rangle
_{Q}\right\vert }.  \label{brutta}
\end{equation}
\end{proposition}

Proof.$\ $We have
\begin{equation*}
e_{0}=\ \underset{\varepsilon \rightarrow 0}{\lim }\ \inf \ \left\{ \frac{%
\frac{1}{2}\left\langle Lu,u\right\rangle _{Q}+K_{Q}(u)}{\left\vert
\left\langle L_{0}u,u\right\rangle _{Q}+K_{0Q}(u)\right\vert }:\text{\ }u\in
X,\ \left\Vert u\right\Vert _{Q}\leq \varepsilon ,\ C_{Q}(u)>0\right\} ;
\end{equation*}%
set $u=\varepsilon v$ with $0<\left\Vert v\right\Vert _{Q}\leq 1;$ then%
\begin{eqnarray*}
e_{0} &=&\ \underset{\varepsilon \rightarrow 0}{\lim }\ \underset{\left\Vert
v\right\Vert _{Q}\leq 1}{\inf }\frac{\frac{1}{2}\left\langle
Lv,v\right\rangle _{Q}+\frac{K_{Q}(\varepsilon v)}{\varepsilon ^{2}}}{%
\left\vert \left\langle L_{0}v,v\right\rangle _{Q}+\frac{K_{0Q}(\varepsilon
v)}{\varepsilon ^{2}}\right\vert } \\
&=&\underset{\left\Vert v\right\Vert _{Q}\leq 1}{\inf }\frac{\frac{1}{2}%
\left\langle Lv,v\right\rangle _{Q}}{\left\vert \left\langle
L_{0}v,v\right\rangle _{Q}\right\vert }=\ \underset{u\neq 0}{\inf }\ \frac{%
\frac{1}{2}\left\langle Lu,u\right\rangle _{Q}}{\left\vert \left\langle
L_{0}u,u\right\rangle _{Q}\right\vert }
\end{eqnarray*}

$\square $\bigskip

\subsection{Proof of the main theorems}

We shall first prove some technical lemmas.

\begin{lemma}
\label{a}Let $Q$ be defined by (\ref{qu}) and $T_{j}q=q+Aj$ $(q\in Q).$ Then%
\begin{equation}
\mathbb{R}^{N}=\dbigcup_{\ j\in \mathbb{Z}^{N}}T_{j}Q.  \label{ricopro}
\end{equation}
\end{lemma}

\textbf{Proof. }Take a generic\textbf{\ }$x\in \mathbb{R}^{N}$ and set $%
y=A^{-1}x.$ We can decompose $y$ as follows: $y=q_{0}+j$ where $j\in \mathbb{%
Z}^{N}$ and $q_{0}\in Q_{0}$ defined by (\ref{qu0}). Then%
\begin{equation*}
x=Ay=Aq_{0}+Aj=q+Aj=T_{j}q
\end{equation*}%
where $q:=Aq_{0}\in Q.$ Since $x$ is generic the lemma is proved.

$\square $

\begin{lemma}
\label{pipo}The map
\begin{equation*}
\sigma \mapsto \Lambda _{\sigma },
\end{equation*}%
where $\Lambda _{\sigma }$ is defined in (\ref{nuovo}), is lower
semicontinuous. Moreover
\begin{equation}
\forall \sigma \in \mathbb{R}^{+},\ \ \left\vert C\left( u\right)
\right\vert =\sigma \Longrightarrow E\left( u\right) \geq \sigma \Lambda
_{\sigma }  \label{fava}
\end{equation}%
and%
\begin{equation}
\underset{\sigma \rightarrow 0}{\lim }\Lambda _{\sigma }\geq e_{0}
\label{uccello}
\end{equation}
\end{lemma}

\textbf{Proof}. The semicontinuity of $\Lambda _{\sigma }$ is an immediate
consequence of the definition. Moreover, by its definition, we have that%
\begin{equation*}
\Lambda _{\sigma }\leq \ \underset{u\in \mathfrak{M}_{\sigma }}{\inf }%
\Lambda \left( u\right) =\ \underset{u\in \mathfrak{M}_{\sigma }}{\inf }%
\frac{E\left( u\right) }{\sigma }
\end{equation*}%
from which (\ref{fava}) follows. Let us prove (\ref{uccello}). We set%
\begin{equation}
e_{\ast }=\ \underset{\varepsilon \rightarrow 0}{\lim }\ \inf \left\{ \frac{%
E(u)}{\left\vert C(u)\right\vert }:\text{\ }u\in X,\ \left\Vert u\right\Vert
\leq \varepsilon ,\ \left\vert C(u)\right\vert >0\right\} .  \label{ul}
\end{equation}%
Let $u_{n}\in X$ be a sequence such that%
\begin{eqnarray*}
\frac{E(u_{n})}{\left\vert C(u_{n})\right\vert } &=&e_{\ast }+o(1) \\
\left\Vert u_{n}\right\Vert &=&o(1).
\end{eqnarray*}%
We can assume, passing eventually to a subsequence, that
\begin{equation*}
C(u_{n})\geq 0.
\end{equation*}%
If such a subsequence does not exist, we have $C(u_{n})\leq 0$ and we argue
in a similar way.

Now take $\varepsilon >0$, then, for $n$ large enough%
\begin{equation*}
\frac{E\left( u_{n}\right) }{C\left( u_{n}\right) }\leq e_{\ast
}+\varepsilon ,
\end{equation*}%
So, if we set
\begin{equation*}
T_{j}Q=\Omega _{j}(j\in \mathbb{Z}^{N}),
\end{equation*}%
by (\ref{ricopro}) (see also section \ref{gf}) we have%
\begin{equation}
e_{\ast }+\varepsilon \geq \frac{E\left( u_{n}\right) }{C\left( u_{n}\right)
}=\frac{\sum_{j\in \mathbb{Z}^{N}}E_{\Omega _{j}}\left( u_{n}\right) }{%
\sum_{j\in \mathbb{Z}^{N}}C_{\Omega _{j}}\left( u_{n}\right) }\geq \frac{%
\sum_{j\in \mathbb{I}}E_{\Omega _{j}}\left( u_{n}\right) }{\sum_{j\in
\mathbb{I}}C_{\Omega _{j}}\left( u_{n}\right) }  \label{mamma}
\end{equation}%
where%
\begin{equation*}
\mathbb{I=}\left\{ j\in \mathbb{Z}^{N}:C_{\Omega _{j}}\left( u_{n}\right)
>0\right\} .
\end{equation*}%
Now, for every $n$ large, it is possible to take $j_{n}\in \mathbb{I}$ such
that%
\begin{equation}
\frac{E_{\Omega _{j_{n}}}\left( u_{n}\right) }{C_{\Omega _{j_{n}}}\left(
u_{n}\right) }\leq e_{\ast }+\varepsilon .  \label{bin}
\end{equation}%
To show this, we argue indirectly and assume that
\begin{equation}
\forall j\in \mathbb{I},\ \frac{E_{\Omega _{j}}\left( u_{n}\right) }{%
C_{\Omega _{j}}\left( u_{n}\right) }>e_{\ast }+\varepsilon ,  \label{bil}
\end{equation}%
then you get
\begin{equation}
\frac{\sum_{j\in \mathbb{I}}E_{\Omega _{j}}\left( u_{n}\right) }{\sum_{j\in
\mathbb{I}}C_{\Omega _{j}}\left( u_{n}\right) }>\frac{\sum_{j\in \mathbb{I}%
}\left( e_{\ast }+\varepsilon \right) C_{\Omega _{j}}\left( u_{n}\right) }{%
\sum_{j\in \mathbb{I}}C_{\Omega _{j}}\left( u_{n}\right) }=e_{\ast
}+\varepsilon .  \label{trin}
\end{equation}%
This contradicts (\ref{mamma}).Now set
\begin{equation*}
v_{n}(x)=u_{n}(x+Aj_{n}).
\end{equation*}%
Then (\ref{bin}) gives%
\begin{equation}
\frac{E_{Q}\left( v_{n}\right) }{C_{Q}\left( v_{n}\right) }\leq e_{\ast
}+\varepsilon .  \label{frim}
\end{equation}%
Since $\left\Vert u_{n}\right\Vert $ and consequently also $\left\Vert
u_{n}\right\Vert _{Q}$ are infinitesimal, from (\ref{frim}) and the
definition of $e_{0}$ we obtain that%
\begin{equation*}
e_{0}\leq \frac{E_{Q}\left( v_{n}\right) }{C_{Q}\left( v_{n}\right) }\leq
e_{\ast }+\varepsilon .
\end{equation*}%
and so%
\begin{equation}
e_{0}\leq e_{\ast }  \label{cazz}
\end{equation}

Now set $L=\ \underset{\sigma \rightarrow 0}{\lim }\Lambda _{\sigma };$ then
there exists a sequence $u_{n}\ $such that $\left\vert C\left( u_{n}\right)
\right\vert \rightarrow 0$ and $\frac{E\left( u_{n}\right) }{\left\vert
C\left( u_{n}\right) \right\vert }\rightarrow L.$ If $E\left( u_{n}\right)
\geq b>0,$ then $L=+\infty $ and (\ref{uccello}) holds. Since
\begin{equation*}
E\left( u_{n}\right) =\frac{1}{2}\left\Vert u_{n}\right\Vert ^{2}+K(u_{n})
\end{equation*}%
if $E\left( u_{n}\right) \rightarrow 0,$ we have that $\left\Vert
u_{n}\right\Vert \ \rightarrow 0$. So, by the definition of $e_{\ast }$ and (%
\ref{cazz}), we have that $L\geq e_{\ast }\geq e_{0}$.

$\square $

\begin{lemma}
\label{sp}\textbf{(Splitting property)} Let $E$ and $C$ be as in (\ref{elle}%
) and (\ref{cu}) and assume that E-3, E-4 are satisfied. Let $%
w_{n}\rightharpoonup 0$ weakly and let $u\in X;$ then%
\begin{eqnarray*}
E\left( u+w_{n}\right) &=&E\left( u\right) +E\left( w_{n}\right) +o(1) \\
C\left( u+w_{n}\right) &=&C\left( u\right) +C\left( w_{n}\right) +o(1)
\end{eqnarray*}
\end{lemma}

\textbf{Proof. }We have to show that $\underset{n\rightarrow \infty }{\lim }%
\left\vert E\left( u+w_{n}\right) -E\left( u\right) -E\left( w_{n}\right)
\right\vert =0.$ By (\ref{elle}), we have that%
\begin{eqnarray*}
&&\underset{n\rightarrow \infty }{\lim }\left\vert E\left( u+w_{n}\right)
-E\left( u\right) -E\left( w_{n}\right) \right\vert \\
&\leq &\ \underset{n\rightarrow \infty }{\lim \frac{1}{2}}\left\vert
\left\Vert u+w_{n}\right\Vert ^{2}-\left\Vert u\right\Vert ^{2}-\left\Vert
w_{n}\right\Vert ^{2}\right\vert +\ \underset{n\rightarrow \infty }{\lim }%
\left\vert K\left( u+w_{n}\right) -K\left( u\right) -K\left( w_{n}\right)
\right\vert .
\end{eqnarray*}

Let us consider each piece independently:%
\begin{equation*}
\underset{n\rightarrow \infty }{\lim }\left\vert \left\Vert
u+w_{n}\right\Vert ^{2}-\left\Vert u\right\Vert ^{2}-\left\Vert
w_{n}\right\Vert ^{2}\right\vert =\ \underset{n\rightarrow \infty }{\lim }%
\left\vert 2\left( u,w_{n}\right) \right\vert =0.
\end{equation*}

Choose $\varepsilon >0$ and $R=R(\varepsilon )>0$ such that
\begin{equation*}
E_{B_{R}^{c}}\left( u\right) <\varepsilon ,\ K_{B_{R}^{c}}\left( u\right)
<\varepsilon \text{ and }\left\Vert u\right\Vert _{B_{R}^{c}}<\varepsilon
\end{equation*}%
$\ $where $\ $%
\begin{equation*}
B_{R}^{c}=\mathbb{R}^{N}-B_{R}(0)\text{ and }B_{R}(0)=\left\{ x\in \mathbb{R}%
^{N}:\left\vert x\right\vert <R\right\} .
\end{equation*}

Then, by the local compactness assumption E-3 (see section \ref{gf}), we
have that%
\begin{eqnarray*}
&&\underset{n\rightarrow \infty }{\lim }\left\vert K\left( u+w_{n}\right)
-K\left( u\right) -K\left( w_{n}\right) \right\vert \\
&=&\ \underset{n\rightarrow \infty }{\lim }\left\vert K_{B_{R}^{c}}\left(
u+w_{n}\right) +K_{B_{R}}\left( u+w_{n}\right) -K_{B_{R}^{c}}\left( u\right)
-K_{B_{R}}\left( u\right) -K_{B_{R}^{c}}\left( w_{n}\right) -K_{B_{R}}\left(
w_{n}\right) \right\vert \\
&\mathbb{=}&\ \underset{n\rightarrow \infty }{\lim }\left\vert
K_{B_{R}^{c}}\left( u+w_{n}\right) -K_{B_{R}^{c}}\left( u\right)
-K_{B_{R}^{c}}\left( w_{n}\right) \right\vert \\
&\leq &\ \underset{n\rightarrow \infty }{\lim }\left\vert
K_{B_{R}^{c}}\left( u+w_{n}\right) -K_{B_{R}^{c}}\left( w_{n}\right)
\right\vert +\varepsilon .
\end{eqnarray*}

By (E-4) and the intermediate value theorem we have that for a suitable $%
\theta \in (0,1)$
\begin{equation*}
\left\vert K_{B_{R}^{c}}\left( u+w_{n}\right) -K_{B_{R}^{c}}\left(
w_{n}\right) \right\vert \leq \left\Vert K_{B_{R}^{c}}^{\prime }\left(
\theta u+\left( 1-\theta \right) w_{n}\right) \right\Vert _{X^{\prime
}\left( B_{R}^{c}\right) }\cdot \left\Vert u\right\Vert _{B_{R}^{c}}\leq
M\cdot \varepsilon
\end{equation*}%
Then%
\begin{equation*}
\underset{n\rightarrow \infty }{\lim }\left\vert K\left( u+w_{n}\right)
-K\left( u\right) -K\left( w_{n}\right) \right\vert \leq \varepsilon +M\cdot
\varepsilon
\end{equation*}%
and since $\varepsilon $ is arbitrary, this limit is 0. Then we have proved
the splitting property for $E.$ The splitting property for $C$ is obtained
arguing in the same way we did with $K.$

$\square $

\begin{lemma}
\label{abstract2}Assume (E-0,..,E-4) and let $\sigma ^{+}$ satisfy the
following assumptions:%
\begin{equation}
\exists \sigma \leq \sigma ^{+}:\ \Lambda _{\sigma }<e_{0}  \label{mona}
\end{equation}%
\begin{equation}
\forall \sigma \geq \sigma ^{+}:\ \sigma ^{+}\Lambda _{_{\sigma ^{+}}}\leq
\sigma \Lambda _{\sigma }.  \label{UH}
\end{equation}%
Then, there exists $\bar{\sigma}\in \left( 0,\sigma ^{+}\right] $ such that
\begin{equation}
\Gamma _{_{\bar{\sigma}}}\neq \varnothing  \label{figata}
\end{equation}%
where $\Gamma _{_{\bar{\sigma}}}$ is as in definition \ref{dss}. Moreover

\begin{itemize}
\item if $u_{n}$ is a sequence such that $\Lambda \left( u_{n}\right)
\rightarrow \Lambda _{_{\bar{\sigma}}}$ and $\left\vert C\left( u_{n}\right)
\right\vert \rightarrow \bar{\sigma}$ then
\begin{equation}
dist(u_{n},\Gamma _{_{\bar{\sigma}}})\rightarrow 0  \label{figata0}
\end{equation}

\item $\Gamma _{_{\bar{\sigma}}}$ has the structure in (\ref{is}), namely%
\begin{equation}
\Gamma _{_{\bar{\sigma}}}=\left\{ u(x-x_{0}):u(x)\in \mathcal{K},\ x_{0}\in
F\right\}  \label{figata1}
\end{equation}%
with $\mathcal{K}$ compact and $F\subset \mathbb{R}^{N}$ is a closed set
such that $F=j+F,\ \forall j\in \mathbb{Z}^{N}.$

\item Any $u\in \Gamma _{_{\bar{\sigma}}}$ solves the eigenvalue problem%
\begin{equation}
E^{\prime }(u)=\lambda C^{\prime }(u).  \label{figata2}
\end{equation}
\end{itemize}
\end{lemma}

\textbf{Proof .} By (\ref{mona}) and Lemma \ref{pipo}, we have that $%
\underset{\sigma \in \left( 0,\sigma ^{+}\right] }{\min \Lambda _{\sigma }}$
exists. Let $\bar{\sigma}\in \left( 0,\sigma ^{+}\right] $ be such that%
\begin{equation}
\Lambda _{_{\bar{\sigma}}}=\ \underset{\sigma \in \left( 0,\sigma ^{+}\right]
}{\min }\Lambda _{\sigma };  \label{topata}
\end{equation}%
by (\ref{mona}), we have that
\begin{equation}
\Lambda _{_{\bar{\sigma}}}<e_{0}.  \label{HH}
\end{equation}

Let $u_{n}\in X$ be a sequence such that
\begin{eqnarray}
\Lambda \left( u_{n}\right) &=&\Lambda _{_{\bar{\sigma}}}+o(1)  \label{A} \\
\left\vert C\left( u_{n}\right) \right\vert &=&\bar{\sigma}+o(1).  \label{B}
\end{eqnarray}

In order to fix the ideas, we may assume that
\begin{equation}
C\left( u_{n}\right) =\bar{\sigma}+o(1).  \label{prima}
\end{equation}

If, on the contrary no subsequence of $C\left( u_{n}\right) $ converge to $%
\bar{\sigma}$, then
\begin{equation*}
C\left( u_{n}\right) \rightarrow -\bar{\sigma}
\end{equation*}
and we argue in a similar way. The proof consists of two steps.

Step1. We prove that for a suitable sequence $\left\{ z_{n}\right\} \subset
\mathbb{Z}^{N}$ we have
\begin{equation*}
u_{n}(x)=\bar{u}(x-Az_{n})+w_{n}(x-Az_{n})
\end{equation*}%
where $\bar{u}\neq 0$ and $w_{n}(x)\rightharpoonup 0$ weakly in $X.\ $

We decompose $\mathbb{R}^{N}$ as in (\ref{ricopro}). Take $\varepsilon >0$,
then, for $n$ large enough%
\begin{equation}
\frac{E\left( u_{n}\right) }{C\left( u_{n}\right) }\leq \Lambda _{_{\bar{%
\sigma}}}+\varepsilon .  \label{base}
\end{equation}

Arguing as in (\ref{mamma}) and (\ref{bin}) (replacing $e_{\ast }$ with $%
\Lambda _{_{\bar{\sigma}}}$), for $n$ large, it is possible to take $%
j_{n}\in \mathbb{I}$ such that
\begin{equation}
\frac{E_{\Omega _{j_{n}}}\left( u_{n}\right) }{C_{\Omega _{j_{n}}}\left(
u_{n}\right) }\leq \Lambda _{_{\bar{\sigma}}}+\varepsilon .  \label{5}
\end{equation}

We set
\begin{equation}
v_{n}(x)=u_{n}(x+Aj_{n}).  \label{suppl}
\end{equation}

By (\ref{A}) and (\ref{B}), $E(u_{n})$ and $C(u_{n})$ are bounded. So, also $%
E(v_{n})$ and $C(v_{n})$ are bounded and, by (E-2), $\left\Vert
v_{n}\right\Vert $ is bounded. Let $\bar{u}$ be the weak limit of $v_{n}.$
We want to show that $\bar{u}\neq 0.$

Clearly $C_{Q}\left( v_{n}\right) =C_{\Omega _{j_{n}}}\left( u_{n}\right) .$
Then, since $j_{n}\in \mathbb{I}$, we have that $C_{Q}\left( v_{n}\right) >0$
and, for $n$ large, by (\ref{5}) we have%
\begin{equation}
\frac{E_{Q}\left( v_{n}\right) }{C_{Q}\left( v_{n}\right) }\leq \Lambda _{_{%
\bar{\sigma}}}+\varepsilon .  \label{bab}
\end{equation}

We claim that the sequence $C_{Q}\left( v_{n}\right) \ $does not converge to
$0$; in fact if $C_{Q}\left( v_{n}\right) \rightarrow 0,$ then, by (\ref{bab}%
), we have that $E_{Q}\left( v_{n}\right) \rightarrow 0.$ Since
\begin{equation*}
E_{Q}\left( v_{n}\right) =\frac{1}{2}\left\Vert v_{n}\right\Vert
_{Q}^{2}+K_{Q}(v_{n}),
\end{equation*}
we have that $\left\Vert v_{n}\right\Vert _{Q}\rightarrow 0;$ so, by
definition of $e_{0}$, and by (\ref{bab}), we have%
\begin{equation*}
\Lambda _{_{\bar{\sigma}}}+\varepsilon \geq \ \underset{n\rightarrow \infty }%
{\lim }\ \frac{E_{Q}\left( v_{n}\right) }{C_{Q}\left( v_{n}\right) }\geq
e_{0}
\end{equation*}%
and this fact contradicts (\ref{HH}) if $\varepsilon >0$ is small enough.

Since $C_{Q}\left( v_{n}\right) \ $does not converge to $0$, by (E-3) with $%
\Omega =Q$, we have that $C_{Q}\left( \bar{u}\right) >0\ $and we can
conclude that $\bar{u}\neq 0.$ Now set
\begin{equation*}
w_{n}=v_{n}-\bar{u}
\end{equation*}%
and so $w_{n}(x)\rightharpoonup 0$ weakly in $X.$

Step 2. Next we will prove that
\begin{equation*}
v_{n}\rightarrow \bar{u}\text{ strongly in }X
\end{equation*}%
namely that $w_{n}\rightarrow 0$ strongly in $X.$ So, by (\ref{elle}), it
will be enough to show that%
\begin{equation}
E\left( w_{n}\right) \rightarrow 0.  \label{e}
\end{equation}

By (\ref{A}), (\ref{B}) and lemma \ref{sp}%
\begin{equation}
\Lambda _{\bar{\sigma}}=\frac{E\left( \bar{u}+w_{n}\right) }{C\left( \bar{u}%
+w_{n}\right) }+o(1)=\frac{E\left( \bar{u}\right) +E\left( w_{n}\right) }{%
\bar{\sigma}}+o\left( 1\right)  \label{primo}
\end{equation}

and so
\begin{equation}
E\left( \bar{u}\right) +E\left( w_{n}\right) =\bar{\sigma}\Lambda _{\bar{%
\sigma}}+o\left( 1\right) .  \label{toga}
\end{equation}

Now we set%
\begin{eqnarray*}
\sigma _{1} &=&\left\vert C\left( \bar{u}\right) \right\vert \\
\sigma _{2} &=&\lim \ \left\vert C\left( w_{n}\right) \right\vert .
\end{eqnarray*}

We consider three cases.

Case 1: $\left\vert C\left( \bar{u}\right) \right\vert =\sigma _{1}\geq
\sigma ^{+}.\ $Then
\begin{eqnarray*}
E\left( \bar{u}\right) &\geq &\sigma _{1}\Lambda _{\sigma _{1}}\ \text{(by (%
\ref{fava}))} \\
&\geq &\sigma ^{+}\Lambda _{\sigma ^{+}}\ \text{(by (\ref{UH}))} \\
&\geq &\sigma ^{+}\Lambda _{\bar{\sigma}}\ \text{(by (\ref{topata}))} \\
&\geq &\bar{\sigma}\Lambda _{\bar{\sigma}}
\end{eqnarray*}%
and by (\ref{toga})%
\begin{equation*}
E\left( w_{n}\right) =\bar{\sigma}\Lambda _{\bar{\sigma}}+o\left( 1\right)
-E\left( \bar{u}\right) \leq o\left( 1\right) ,
\end{equation*}%
and so $E\left( w_{n}\right) \rightarrow 0.$

Case 2: $\sigma _{2}=\left\vert C\left( w_{n}\right) +o(1)\right\vert \geq
\sigma ^{+}.$ Then%
\begin{eqnarray*}
E\left( w_{n}\right) &\geq &\left\vert C\left( w_{n}\right) \right\vert
\Lambda _{\left\vert C\left( w_{n}\right) \right\vert }\ \text{(by (\ref%
{fava}))} \\
&\geq &\sigma _{2}\Lambda _{\sigma _{2}}+o(1)\ \text{(by lemma \ref{pipo})}
\\
&\geq &\sigma ^{+}\Lambda _{\sigma ^{+}}+o(1)\ \text{(by (\ref{UH}))} \\
&\geq &\sigma ^{+}\Lambda _{\bar{\sigma}}+o(1)\ \text{(by (\ref{topata}))} \\
&\geq &\bar{\sigma}\Lambda _{\bar{\sigma}}+o(1).
\end{eqnarray*}%
Then by (\ref{toga}) you get%
\begin{equation*}
\bar{\sigma}\Lambda _{\bar{\sigma}}=E\left( \bar{u}\right) +E\left(
w_{n}\right) +o\left( 1\right) \geq E\left( \bar{u}\right) +\bar{\sigma}%
\Lambda _{\bar{\sigma}}+o(1)
\end{equation*}%
and this is a contradiction since $E\left( \bar{u}\right) >0;$ thus case 2
cannot occur.

Case 3: $\sigma _{1},\sigma _{2}\leq \sigma ^{+}.$ In this case, we have by (%
\ref{toga}) and (\ref{fava})%
\begin{eqnarray*}
\bar{\sigma}\Lambda _{\bar{\sigma}} &=&E\left( \bar{u}\right) +E\left(
w_{n}\right) +o\left( 1\right) \geq \sigma _{1}\Lambda _{\sigma _{1}}+\sigma
_{2}\Lambda _{\sigma _{2}}+o(1) \\
&\geq &\left( \sigma _{1}+\sigma _{2}\right) \Lambda _{\bar{\sigma}}.
\end{eqnarray*}%
Then%
\begin{equation}
\sigma _{1}+\sigma _{2}\leq \bar{\sigma}.  \label{ma}
\end{equation}

Now the opposite inequality can be obtained by splitting the charge as in
lemma \ref{sp} :%
\begin{equation}
\bar{\sigma}=\left\vert C(\bar{u}+w_{n})\right\vert +o(1)\leq \left\vert C(%
\bar{u})\right\vert +\left\vert C(w_{n})\right\vert +o(1)=\sigma _{1}+\sigma
_{2}+o(1).  \label{mi}
\end{equation}%
From (\ref{ma}) and (\ref{mi}) we get
\begin{equation}
\sigma _{1}+\sigma _{2}=\bar{\sigma}.  \label{per}
\end{equation}

Now we claim that
\begin{equation}
\sigma _{1}>0.  \label{claim}
\end{equation}

Arguing by contradiction assume $\sigma _{1}=0,$ then by ($\ref{per}$) we
have $\sigma _{2}=\bar{\sigma}$ and by (\ref{toga}) and (\ref{fava})%
\begin{equation*}
\bar{\sigma}\Lambda _{\bar{\sigma}}=E\left( \bar{u}\right) +E\left(
w_{n}\right) +o\left( 1\right) \geq E\left( \bar{u}\right) +\sigma
_{2}\Lambda _{\sigma _{2}}+o\left( 1\right) =E\left( \bar{u}\right) +\bar{%
\sigma}\Lambda _{\bar{\sigma}}+o\left( 1\right)
\end{equation*}%
and this contradicts $E\left( \bar{u}\right) >0.$

Now it is not restrictive to suppose that
\begin{equation}
\bar{\sigma}=\min \left\{ \sigma :\Lambda _{\sigma }=\underset{\tau \in
\left( 0,\sigma ^{+}\right] }{\min }\Lambda _{\tau }\right\} .
\label{lillina}
\end{equation}%
We claim that $\sigma _{2}=0.$ In fact, arguing by contradiction assume that
$\sigma _{2}>0,\ $then, by ($\ref{per}$), $\sigma _{1}<\bar{\sigma}$ and, by
(\ref{lillina}), $\Lambda _{\sigma _{1}}-\Lambda _{\bar{\sigma}}=\delta >0.$
So we have%
\begin{eqnarray*}
\bar{\sigma}\Lambda _{\bar{\sigma}} &=&E\left( \bar{u}\right) +E\left(
w_{n}\right) +o\left( 1\right) \geq \sigma _{1}\Lambda _{\sigma _{1}}+\sigma
_{2}\Lambda _{\sigma _{2}}+o(1) \\
&\geq &\sigma _{1}\left( \Lambda _{\bar{\sigma}}+\delta \right) +\sigma
_{2}\Lambda _{\bar{\sigma}}+o(1)=\bar{\sigma}\Lambda _{\bar{\sigma}}+\sigma
_{1}\delta +o(1)
\end{eqnarray*}%
and this is a contradiction since $\sigma _{1}\delta >0$, so we have $\sigma
_{2}=0.$

Since $\sigma _{2}=0,$ then $\sigma _{1}=\bar{\sigma},$ and by (\ref{toga})
and (\ref{fava})
\begin{equation*}
E\left( w_{n}\right) =\bar{\sigma}\Lambda _{\bar{\sigma}}-E\left( \bar{u}%
\right) +o\left( 1\right) \leq \bar{\sigma}\Lambda _{\bar{\sigma}}-\sigma
_{1}\Lambda _{\sigma _{1}}+o\left( 1\right) =o\left( 1\right)
\end{equation*}%
from which we get (\ref{e}).

By the preceding results we easily get the conclusions (\ref{figata},...,\ref%
{figata2}). In fact:

-Consider the sequence $v_{n}$ defined in (\ref{suppl}). We have seen in
steps 1, 2 that $v_{n}\rightarrow \bar{u}$ strongly in $X.$ Then, since $E$
and $C$ are continuous, we have
\begin{equation}
\frac{E(v_{n})}{C(v_{n})}=\frac{E(\bar{u})}{C(\bar{u})}+o(1)=\frac{E(\bar{u})%
}{\bar{\sigma}}+o(1).  \label{na}
\end{equation}

Moreover by (\ref{A})
\begin{equation}
\frac{E(v_{n})}{C(v_{n})}=\Lambda _{\bar{\sigma}}+o(1)\leq \inf \left\{
\frac{E(u)}{\bar{\sigma}}:C(u)=\bar{\sigma}\right\} +o(1).  \label{ne}
\end{equation}

From (\ref{na}) and (\ref{ne}) we deduce that
\begin{equation}
\bar{u}\in \Gamma _{\bar{\sigma}}.  \label{fig}
\end{equation}

-By steps 1, 2 and (\ref{fig}), we clearly get (\ref{figata0}). Moreover, if
we take a sequence $\left\{ u_{n}\right\} \subset \Gamma _{\bar{\sigma}}$,
by using again steps 1,2, we get that there exists a subsequence, which we
continue to call $u_{n},$ and $\left\{ j_{n}\right\} $ $\subset F$ such that
\begin{equation*}
v_{n}\rightarrow \bar{u}\in \Gamma _{\bar{\sigma}}\text{ strongly in }X,%
\text{ }v_{n}(x)=u_{n}(Aj_{n}+x).
\end{equation*}%
Then also (\ref{figata1}) holds.

Finally (\ref{figata2}) clearly follows by the definition of $\Gamma _{\bar{%
\sigma}}$ .

$\square $

\textbf{Proof of Th. \ref{due}.} We prove that the assumptions (\ref{mona})
and (\ref{UH}) of Lemma \ref{abstract2} are satisfied.

First, we observe that, by (\ref{no}), $\Lambda _{\sigma }\geq \Lambda
^{\ast }>0,$ then%
\begin{equation}
\sigma \Lambda _{\sigma }\rightarrow \infty \text{ for }\sigma \rightarrow
\infty .  \label{beg}
\end{equation}

Now set%
\begin{equation}
\tau _{n}=\sup \left\{ \sigma :\sigma \Lambda _{\sigma }\geq n\right\} .
\label{numero}
\end{equation}%
Then, by definition%
\begin{equation}
\tau _{n}\Lambda _{\tau _{n}}\leq n  \label{nume}
\end{equation}

By (\ref{beg})%
\begin{equation}
\tau _{n}\in \mathbb{R}\text{ and }\tau _{n}\rightarrow \infty \text{ for }%
n\rightarrow \infty .  \label{big}
\end{equation}

Now by (\ref{no}) there exists $u_{0}\in X$ such that%
\begin{equation}
\Lambda (u_{0})<e_{0}.  \label{us}
\end{equation}

By (\ref{big}) there exists $\bar{n}$ such that
\begin{equation*}
\tau _{\bar{n}}\geq \left\vert C\left( u_{0}\right) \right\vert
\end{equation*}%
and by (\ref{numero}) and (\ref{nume})
\begin{equation}
\text{ }\tau _{\bar{n}}\Lambda _{\tau _{\bar{n}}}\leq \bar{n}\leq \sigma
\Lambda _{\sigma }\text{ for }\sigma \geq \tau _{\bar{n}}.  \label{his}
\end{equation}

Set $\sigma _{+}=\tau _{\bar{n}}$, then by (\ref{us}) and (\ref{his}) we get
\begin{eqnarray*}
\sigma &=&\left\vert C\left( u_{0}\right) \right\vert \leq \sigma ^{+},\text{
}\Lambda _{\sigma }<e_{0} \\
\text{ }\sigma ^{+}\Lambda _{\sigma ^{+}} &\leq &\sigma \Lambda _{\sigma }%
\text{ for }\sigma \geq \sigma ^{+}.
\end{eqnarray*}%
Then the assumptions (\ref{mona}) and (\ref{UH}) of Lemma (\ref{abstract2})
are satisfied.

$\square $

\bigskip

\begin{remark}
\label{figo}By the proof of this theorem, we can see that the assumption $%
\Lambda ^{\ast }>0$ is used only to get (\ref{beg}). This assumption can be
replaced by the following one%
\begin{equation}
\left\Vert u_{n}\right\Vert \rightarrow \infty \Rightarrow E\left(
u_{n}\right) \rightarrow \infty  \label{ei}
\end{equation}%
In fact (\ref{ei}) implies (\ref{beg}). To show this, we argue indirectly
and assume that there exists a sequence $\sigma _{n}\rightarrow \infty $
such that $\sigma _{n}\Lambda _{\sigma _{n}}$ is bounded; so there exists a
sequence $u_{n}$ such that
\begin{equation}
\left\vert C\left( u_{n}\right) \right\vert \rightarrow \infty
\label{racchia}
\end{equation}%
and
\begin{equation}
\left\vert C\left( u_{n}\right) \right\vert \Lambda \left( u_{n}\right)
=E\left( u_{n}\right) \text{ is bounded.}  \label{bono}
\end{equation}%
By (\ref{racchia}) and (E-4)$,$ we have that (for a subsequence) $\left\Vert
u_{n}\right\Vert \rightarrow \infty ;$ then, by (\ref{ei}), $E\left(
u_{n}\right) \rightarrow \infty .$ This contradicts (\ref{bono}), then we
conclude that (\ref{beg}) holds.
\end{remark}

\textbf{Proof of Th. \ref{tre}.} Arguing as in the proof of Th. \ref{due},
there exists $\sigma _{0}^{+}>0$ such that
\begin{equation}
\Lambda \left( \bar{u}\right) =\ \underset{\left\vert C\left( u\right)
\right\vert \leq \sigma _{0}^{+}}{\min }\Lambda \left( u\right) .
\label{lume}
\end{equation}%
for a suitable $\bar{u}.$

We shall show that%
\begin{equation}
\left\vert C\left( \bar{u}\right) \right\vert =\sigma _{0}^{+}.  \label{pelo}
\end{equation}

Arguing by contradiction, assume that $\left\vert C\left( \bar{u}\right)
\right\vert <\sigma _{0}^{+}.$ Then, since $C(R_{\theta }u)\ $and $\Lambda
\left( R_{\theta }u\right) $ are respectively increasing and decreasing in $%
\theta $, for $\varepsilon >0$ small enough and $1<\vartheta <1+\varepsilon $%
, we have
\begin{equation}
\sigma _{0}^{+}\geq \left\vert C\left( R_{\theta }\bar{u}\right) \right\vert
>\left\vert C\left( \bar{u}\right) \right\vert  \label{kappa}
\end{equation}

\begin{equation}
\Lambda \left( R_{\theta }\bar{u}\right) <\Lambda \left( \bar{u}\right) .
\label{kappone}
\end{equation}%
Clearly (\ref{kappa}) and (\ref{kappone}) contradict (\ref{lume}). So (\ref%
{pelo}) holds.

\ Now set $\sigma _{0}=\sigma _{0}^{+}$ and take any other $\sigma ^{+}\geq
\sigma _{0}.$ Clearly (\ref{mona}) and (\ref{UH}) hold and we can argue as
before.

$\square $

\bigskip

\subsection{Dynamical consequences of the main theorem\label{cmt}}

The above theorems can be applied to the case in which $\left( X,\left\Vert
\cdot \right\Vert \right) $ is the state space of a dynamical system $(X,T)$
and it proves the existence of hylomorphic solitons; more exactly we have:

\begin{theorem}
Let $(X,T)$ be a dynamical system and let $E$ and $C$ be the energy and the
charge. If $X,E$ and $C$ are as in section \ref{gf} and satisfy the
assumptions of theorem \ref{due}, then $(X,T)$ has hylomorphic solitons.
Moreover, if also the assumptions of Th. \ref{tre} are satisfied, there
exists $\sigma _{0}$ such that there are solitons for any charge $\bar{\sigma%
}\geq \sigma _{0}.$
\end{theorem}

\textbf{Proof.} We consider Def. \ref{dss}. We set
\begin{equation*}
\Gamma _{\sigma }=\left\{ u\in \mathfrak{M}_{\sigma }:E(u)=c_{\sigma
}\right\}
\end{equation*}%
with%
\begin{equation*}
c_{\sigma }=\ \underset{u\in \mathfrak{M}_{\sigma }}{\min }E(u).
\end{equation*}%
By theorem \ref{due} $\Gamma _{\sigma }\neq \varnothing .$ In order to prove
the existence of solitons we need to prove (ii) and (iii) of definition \ref%
{ds}. (ii) follows by (\ref{figata1}).

In order to prove stability, we use the Lyapunov criterium; we define the
Lyapunov function $V:X\rightarrow \mathbb{R}$ as follows
\begin{equation*}
V(u):=\left( E(u)-c_{\sigma }\right) ^{2}+\left( \left\vert C(u)\right\vert
-\sigma \right) ^{2};
\end{equation*}%
then by (\ref{figata0})
\begin{equation}
V(u_{n})\rightarrow 0\Longrightarrow d\left( u_{n},\Gamma \right)
\rightarrow 0.  \label{vto0}
\end{equation}%
Then, by the Lyapunov stability theorem $\Gamma $ is stable.

The second statements follows directly from Th. \ref{tre}.

$\square $

\section{The nonlinear Schr\"{o}dinger equation\label{SSE}}

We are interested to the nonlinear Schr\"{o}dinger equation:
\begin{equation}
i\frac{\partial \psi }{\partial t}=-\frac{1}{2}\Delta \psi +V(x)\psi +\frac{1%
}{2}W^{\prime }(\psi )  \label{NSV}
\end{equation}%
where $\psi :\mathbb{R}^{N}\mathbb{\rightarrow C}$ $(N\geq 3),\ V:\mathbb{R}%
^{N}\mathbb{\rightarrow R},\ W:\mathbb{C\rightarrow R}$ such that $%
W(s)=F(\left\vert s\right\vert )$ for some smooth function $F:\left[
0,\infty \right) \rightarrow \mathbb{R}$ and
\begin{equation}
W^{\prime }(s)=\frac{\partial W}{\partial s_{1}}+i\frac{\partial W}{\partial
s_{2}},\ \ \ s=s_{1}+is_{2}  \label{w'}
\end{equation}%
namely
\begin{equation*}
W^{\prime }(s)=F^{\prime }(\left\vert s\right\vert )\frac{s}{\left\vert
s\right\vert }.
\end{equation*}%
.

Equation (\ref{NSV}) is the Euler-Lagrange equation relative to the
Lagrangian density
\begin{equation}
\mathcal{L}=\func{Re}\left( i\partial _{t}\psi \overline{\psi }\right) -%
\frac{1}{2}\left\vert \nabla \psi \right\vert ^{2}-V(x)\left\vert \psi
\right\vert ^{2}-W\left( \psi \right)  \label{lagr}
\end{equation}

\subsection{Existence results}

We assume that $W$ has the following form%
\begin{equation}
W(s)=\frac{1}{2}\ h^{2}s^{2}+N(s)  \label{NO}
\end{equation}%
where $h^{2}=W^{\prime \prime }(0)$ and $N(s)=o(s^{2}).$ We make the
following assumptions on $W$:

\begin{itemize}
\item (W-i) \textbf{(Positivity}) $W(s)\geq 0$

\item (W-ii) \textbf{(Nondegeneracy}) $W=$ $W(s)$ ( $s\geq 0)$ is $C^{2}$
near the origin with $W(0)=W^{\prime }(0)=0;\;W^{\prime \prime }(0)>0$

\item (W-iii) \textbf{(Hylomorphy}) $0<\inf \frac{W(s)}{\frac{1}{2}s^{2}}%
<W^{\prime \prime }(0)$\

\item (W-iiii) \textbf{(Growth condition}) there are constants $%
c_{1},c_{2}>0,$ $2<p<2N/(N-2)$ such that for any $s>0:$%
\begin{equation*}
|N^{\prime }(s)|\ \leq c_{1}s^{p-1}+c_{2}s^{2-\frac{2}{p}}.
\end{equation*}
\end{itemize}

If we set%
\begin{equation}
\alpha ^{2}=\inf \frac{W(s)}{\frac{1}{2}s^{2}},  \label{alpha}
\end{equation}%
then the hylomorphy assumption (W-iii) reads
\begin{equation}
0<\alpha <h.  \label{illo}
\end{equation}

This assumption implies that%
\begin{equation*}
\exists \bar{s}:N(\bar{s})<0.
\end{equation*}%
We make the following assumptions on $V:$

\begin{itemize}
\item (V-i) \textbf{(Positivity}) $V\geq 0$ and $V\in L^{\infty .}$

\item (V-ii) \textbf{(Lattice invariance}) There exists an $N\times N$
invertible matrix $A$ such that
\begin{equation}
V(x)=V(x-Az)  \label{vii}
\end{equation}
for all $x\in \mathbb{R}^{N}$ and $z\in \mathbb{Z}^{N}.$
\end{itemize}

Here we want to use the results of the previous sections to study (\ref{NSV}%
). In this case the state $u$ coincides with $\psi $ and the general
framework of the previous sections takes the following form:

\begin{equation*}
X=H^{1}(\mathbb{R}^{N})
\end{equation*}%
where $H^{1}(\mathbb{R}^{N})$ is the usual Sobolev space and
\begin{eqnarray}
E\left( u\right) &=&\int \left( \frac{1}{2}\left\vert \nabla u\right\vert
^{2}+V(x)\left\vert u\right\vert ^{2}+W\left( u\right) \right) dx  \label{en}
\\
&=&\int \left( \frac{1}{2}\left\vert \nabla u\right\vert ^{2}+\frac{%
h^{2}u^{2}}{2}+V(x)\left\vert u\right\vert ^{2}\right) dx+\int N(u)dx;
\end{eqnarray}

\begin{equation}
C\left( u\right) =\int u^{2}dx  \label{cn}
\end{equation}%
\begin{equation*}
\left\Vert u\right\Vert ^{2}=\int \left( \left\vert \nabla u\right\vert
^{2}+au^{2}\right) dx.
\end{equation*}%
Then the energy $E$ and the hylenic charge $C$ have the form (\ref{elle})
and (\ref{cu}) respectively. We shall prove the following theorem

\begin{theorem}
\label{sch} Assume that $W$ satisfies W-i),...W-iiii) and that $V$ satisfies
V-i), V-ii)$.$ Moreover assume that
\begin{equation}
\frac{\alpha ^{2}}{2}+\left\Vert V\right\Vert _{L^{\infty }}<\frac{h^{2}}{2}
\label{cos}
\end{equation}%
where $\alpha $ and $h$ have been introduced in (\ref{alpha}) and (\ref{NO}%
). Then equation (\ref{NSV}) admits hylomorphic solitons (see definition \ref%
{dss}).
\end{theorem}

\begin{remark}
Observe that, when $V=0,$ assumption (\ref{cos}) reduces to the request $%
\alpha <h,$ which is the "usual" hylomorphy condition (see \cite{BBBM}, \cite%
{hylo}, \cite{befo08}, \cite{milano}). Moreover, in this case it is possible
to apply Th. \ref{tre} and to get the existence of solitons for any
sufficiently large charge.
\end{remark}

\begin{remark}
Actually, the assumptions (W-i,...,W-iiii) are not the most general. For
example the positivity assumption is not necessary. In the case $V=0,$ we
refer to \cite{BBGM}. If $V\neq 0$, we do not know whether the assumptions
used in \cite{BBGM} are sufficient.
\end{remark}

We first obtain some estimates on $e_{0}$ and $\Lambda _{\ast }$ defined by (%
\ref{brutta}) and (\ref{brutta1}).

\begin{lemma}
\label{comp}Assume that $W$ satisfies (W-i,...W-iiii) and that $V$ satisfies
(V-i, V-ii). Then
\begin{eqnarray}
\frac{h^{2}}{2} &\leq &e_{0}\leq \frac{h^{2}}{2}+\left\Vert V\right\Vert
_{L^{\infty }}  \label{lula} \\
\frac{\alpha ^{2}}{2} &\leq &\Lambda _{\ast }\leq \frac{\alpha ^{2}}{2}%
+\left\Vert V\right\Vert _{L^{\infty }}  \label{lulaa}
\end{eqnarray}
\end{lemma}

\textbf{Proof}. By using (\ref{brutta}), we clearly deduce that (\ref{lula})
holds.

Now we prove (\ref{lulaa}). First we show that:%
\begin{equation}
\Lambda _{\ast }\geq \frac{\alpha ^{2}}{2}.  \label{ba}
\end{equation}

In fact, by using (\ref{alpha}), we get

\begin{align*}
\Lambda _{\ast }& =\underset{u}{\ \inf }\ \frac{E\left( u\right) }{%
\left\vert C\left( u\right) \right\vert }=\underset{u}{\ \inf }\ \frac{\int
\left( \frac{1}{2}\left\vert \nabla u\right\vert ^{2}+V(x)\left\vert
u\right\vert ^{2}+W(u)\right) dx}{\int u^{2}dx} \\
& \geq \underset{u}{\ \inf }\frac{\int W(u)dx}{\int u^{2}dx}\geq \underset{u}%
{\ \inf }\frac{\int \frac{1}{2}\alpha ^{2}u^{2}dx}{\int u^{2}dx}=\frac{1}{2}%
\alpha ^{2}.
\end{align*}%
Now we prove that%
\begin{equation}
\Lambda _{\ast }\leq \frac{\alpha ^{2}}{2}+\left\Vert V\right\Vert
_{L^{\infty }}.  \label{baa}
\end{equation}

Take $\varepsilon >0,$ then by (\ref{alpha}), there exists $s_{\varepsilon
}>0$ such that
\begin{equation}
W(s_{\varepsilon })<\frac{1}{2}s_{\varepsilon }^{2}(\alpha ^{2}+\varepsilon
).  \label{inf}
\end{equation}

Let $R>0$ and set
\begin{equation}
u_{\varepsilon ,R}=\left\{
\begin{array}{cc}
s_{\varepsilon } & if\;\;|x|<R \\
0 & if\;\;|x|>R+1 \\
\frac{|x|}{R}s_{\varepsilon }-(\left\vert x\right\vert -R)\frac{R+1}{R}%
s_{\varepsilon } & if\;\;R<|x|<R+1%
\end{array}%
\right.  \label{copia}
\end{equation}

Clearly%
\begin{equation}
\frac{\int \frac{1}{2}\left\vert \nabla u_{\varepsilon ,R}\right\vert ^{2}}{%
\int \left\vert u_{\varepsilon ,R}\right\vert ^{2}dx}\leq O\left( \frac{1}{R}%
\right) .  \label{pre}
\end{equation}%
Then, by (\ref{inf}) and (\ref{pre}) we get

\begin{eqnarray*}
\Lambda _{\ast } &\leq &\frac{\int \left( \frac{1}{2}\left\vert \nabla
u_{\varepsilon ,R}\right\vert ^{2}+V(x)\left\vert u\right\vert
^{2}+W(u_{\varepsilon ,R})\right) dx}{\int \left\vert u_{\varepsilon
,R}\right\vert ^{2}dx} \\
&\leq &\frac{\int_{\left\vert x\right\vert <R}\left( W(u_{\varepsilon
,R})+V(x)\left\vert u_{\varepsilon ,R}\right\vert ^{2}\right) dx}{%
\int_{\left\vert x\right\vert <R}\left\vert u_{\varepsilon ,R}\right\vert
^{2}dx} \\
&&+\frac{\int_{R<\left\vert x\right\vert <R+1}\left( \frac{1}{2}\left\vert
\nabla u_{\varepsilon ,R}\right\vert ^{2}+W(u_{\varepsilon
,R})+V(x)\left\vert u_{\varepsilon ,R}\right\vert ^{2}\right) dx}{%
\int_{\left\vert x\right\vert <R}\left\vert u_{\varepsilon ,R}\right\vert
^{2}dx} \\
&\leq &\frac{\int_{\left\vert x\right\vert <R}\left( W(u_{\varepsilon
,R})+V(x)\left\vert u_{\varepsilon ,R}\right\vert ^{2}\right) dx}{%
\int_{\left\vert x\right\vert <R}\left\vert u_{\varepsilon ,R}\right\vert
^{2}dx}+\frac{c_{1}R^{N-1}}{c_{2}R^{N}} \\
&=&\frac{\int_{\left\vert x\right\vert <R}\left( W(s_{\varepsilon
})+V(x)\left\vert s_{\varepsilon }\right\vert ^{2}\right) dx}{%
\int_{\left\vert x\right\vert <R}\left\vert s_{\varepsilon }\right\vert
^{2}dx}+O\left( \frac{1}{R}\right) \\
&\leq &\frac{\frac{1}{2}s_{\varepsilon }^{2}(\alpha ^{2}+\varepsilon )R^{N}}{%
s_{\varepsilon }^{2}R^{N}}+\frac{\left\Vert V\right\Vert _{L^{\infty
}}s_{\varepsilon }^{2}R^{N}}{s_{\varepsilon }^{2}R^{N}}+O\left( \frac{1}{R}%
\right) =\frac{1}{2}(\alpha ^{2}+\varepsilon )+\left\Vert V\right\Vert
_{L^{\infty }}+O\left( \frac{1}{R}\right)
\end{eqnarray*}

Then, since $\varepsilon >0$ is arbitrary, we easily get (\ref{baa}).
Finally (\ref{lulaa}) follows from (\ref{ba}) and (\ref{baa}).

$\square $

Proof of Theorem \ref{sch}:

By (\ref{lula}), (\ref{lulaa}) and (\ref{cos}) we deduce that $0<\Lambda
_{\ast }<e_{0}.$ It can be shown, by standard calculations (see e.g. \cite%
{befo08}), that under the assumptions W-i),...,W-iiii) and V-i), V-ii), the
functionals $E$ and $C,$ defined by (\ref{en}) and (\ref{cn}), satisfy
(E-0,..,E-4) of section \ref{gf}. Then, by using Theorem \ref{due}, we
deduce that equation (\ref{NSV}) admits hylomorphic solitons. Since these
solitons $u_{0}$ are minimizers of the energy $E$ on the manifold $\left\{
u\in H^{1}(\mathbb{R}^{N}):C(u)=\int u^{2}dx=\sigma \right\} ,$ we get
\begin{equation*}
E^{\prime }(u_{0})=-\omega C^{\prime }(u_{0})
\end{equation*}%
where $\omega $ is a Lagrange multiplier. Then it can be easily seen that $%
u_{0}$ solves (\ref{NSV}) and $u_{0}=\psi _{0}(x)e^{-i\omega t},$ where $%
\omega \in \mathbb{R}$ and $\psi _{0}(x)$ solve the equation
\begin{equation*}
-\frac{1}{2}\Delta \psi _{0}+V(x)\psi _{0}+\frac{1}{2}W^{\prime }(\psi
_{0})=\omega \psi _{0}
\end{equation*}

$\square$

\section{The nonlinear Klein-Gordon equation\label{SKG}}

In this section we will apply th. \ref{due} to the existence of hylomorphic
solitons of the nonlinear Klein-Gordon equation. We point out that the
existence of such solitons for this equation has been recently stated in
\cite{BBBM}. Here we consider the case in which $W$ depends on $x$ and it
has a lattice symmetry.

More exactly, we consider the equation%
\begin{equation}
\square \psi +W^{\prime }(x,\psi )=0  \label{KG}
\end{equation}%
where $\square =\partial _{t}^{2}-\nabla ^{2}$,$\;\psi :\mathbb{R}%
^{N}\rightarrow \mathbb{C}$ ($N\geq 3$) , $W:\mathbb{R}^{N}\times \mathbb{C}%
\rightarrow \mathbb{R}$ and $W^{\prime }$ is the derivative with respect to
the second variable as in (\ref{w'}). We can write $W$ as follows%
\begin{equation}
W(x,s)=\frac{1}{2}\ h(x)^{2}s^{2}+N(x,s),\ \ x\in \mathbb{R}^{N},\ s\in
\mathbb{R}^{+},\ h(x)\in L^{\infty }  \label{h0}
\end{equation}%
where
\begin{equation}
h(x)\geq h_{0}>0  \label{h1}
\end{equation}%
and$\ N(x,s)=o(s^{2}).$ We make the following assumptions on $W$:

\begin{itemize}
\item (NKG-i) \textbf{(Positivity}) $W(x,s)\geq 0$

\item (NKG-ii) \textbf{(Lattice invariance}) There exists an $N\times N$
invertible matrix $A$ such that
\begin{equation}
W(x,s)=W(x-Az,s)  \label{gii}
\end{equation}
for all $x\in \mathbb{R}^{N}$ and $z\in \mathbb{Z}^{N}.$

\item (NKG-iii) \textbf{(Hylomorphy}) $\exists \alpha ,\bar{s}\in \mathbb{R}%
^{+}\mathbb{\ }$such that $W(x,\bar{s})\leq \frac{1}{2}\alpha ^{2}\bar{s}%
^{2} $

\item (NKG-iiii)\textbf{(Growth condition}) there are constants $%
c_{1},c_{2}>0,$ $2<p<2N/(N-2)$ such that for any $s>0:$%
\begin{equation*}
|N^{\prime }(x,s)|\ \leq c_{1}s^{p-1}+c_{2}s^{2-\frac{2}{p}}.
\end{equation*}
\end{itemize}

We shall assume that the initial value problem is well posed for (NKG).

Eq. (\ref{KG}) is the Euler-Lagrange equation of the action functional
\begin{equation}
S(\psi )=\int \left( \frac{1}{2}\left\vert \partial _{t}\psi \right\vert
^{2}-\frac{1}{2}|\nabla \psi |^{2}-W(x,\psi )\right) dxdt.  \label{az}
\end{equation}

The energy and the charge take the following form:
\begin{equation}
E(\psi )=\int \left[ \frac{1}{2}\left\vert \partial _{t}\psi \right\vert
^{2}+\frac{1}{2}\left\vert \nabla \psi \right\vert ^{2}+W(x,\psi )\right] dx
\label{energy}
\end{equation}%
\begin{equation}
C(\psi )=-\func{Re}\int i\partial _{t}\psi \overline{\psi }\;dx.  \label{im1}
\end{equation}

(the sign "minus"in front of the integral is a useful convention).

\subsection{The NKG as a dynamical system}

We set%
\begin{equation*}
X=H^{1}(\mathbb{R}^{N},\mathbb{C})\times L^{2}(\mathbb{R}^{N},\mathbb{C})
\end{equation*}%
and we will denote the generic element of $X$ by $u=(\psi \left( x\right) ,%
\hat{\psi}\left( x\right) );$ then, by the well posedness assumption, for
every $u\in X,$ there is a unique solution $\psi (t,x)$ of (\ref{KG}) such
that%
\begin{eqnarray*}
\psi (0,x) &=&\psi \left( x\right) \\
\partial _{t}\psi (0,x) &=&\hat{\psi}\left( x\right) .
\end{eqnarray*}%
Thus, using our notation, we can write
\begin{equation*}
T_{t}u=U(t,x)=(\psi \left( t,x\right) ,\hat{\psi}\left( t,x\right) )\in
C^{1}(\mathbb{R},X).
\end{equation*}

Using this notation, we can write equation (\ref{KG}) in Hamiltonian form:
\begin{eqnarray}
\partial _{t}\psi &=&\hat{\psi}  \label{ham1} \\
\partial _{t}\hat{\psi} &=&\Delta \psi -W^{\prime }(x,\psi ).  \label{ham2}
\end{eqnarray}%
The energy and the charge, as functionals defined in $X,$ become

\begin{equation}
E(u)=\int \left[ \frac{1}{2}\left\vert \hat{\psi}\right\vert ^{2}+\frac{1}{2}%
\left\vert \nabla \psi \right\vert ^{2}+W(x,\psi )\right] dx
\end{equation}%
\begin{equation}
C(u)=-\func{Re}\int i\hat{\psi}\overline{\psi }\;dx.
\end{equation}%
We shall tacitely assume that $W$ is such that $E,$ $C$ are $C^{1}$ in $X.$

\begin{proposition}
\label{pipi}Let $u_{0}(x)=(\psi _{0}(x),\hat{\psi}_{0}(x))$ $\in X$ be a
critical point of $E$ constrained on the manifold $\mathfrak{M}_{\sigma
}=\left\{ u\in X:C(u)=\sigma \right\} $. Then there exists $\omega \in
\mathbb{R}$ such that $\psi _{0}$ satisfies the equation
\begin{equation}
-\Delta \psi _{0}+W^{\prime }(x,\psi _{0})=\omega ^{2}\psi _{0}  \label{ella}
\end{equation}%
and
\begin{equation}
U(t,x)=\left[
\begin{array}{c}
\psi _{0}(x)e^{-i\omega t} \\
-i\omega \psi _{0}(x)e^{-i\omega t}%
\end{array}%
\right]  \label{evol}
\end{equation}%
solves (\ref{ham1}), (\ref{ham2}).
\end{proposition}

Proof. Clearly
\begin{equation}
E^{\prime }(u_{0})=-\omega C^{\prime }(u_{0})  \label{become}
\end{equation}

where $-\omega $ is a Lagrange multiplier. We now compute the derivatives $%
E^{\prime }(u_{0}),C^{\prime }(u_{0}).$

For all $(v_{0},v_{1})\in X=H^{1}(\mathbb{R}^{N},\mathbb{C})\times L^{2}(%
\mathbb{R}^{N},\mathbb{C}),$ we have%
\begin{equation*}
E^{\prime }(u_{0})\left[
\begin{array}{c}
v_{0} \\
v_{1}%
\end{array}%
\right] =\func{Re}\int \left[ \hat{\psi}_{0}\overline{v_{1}}+\nabla \psi _{0}%
\overline{\nabla v_{0}}+W^{\prime }(x,\psi _{0})\overline{v_{0}}\right] dx
\end{equation*}%
\begin{eqnarray*}
C^{\prime }(u_{0})\left[
\begin{array}{c}
v_{0} \\
v_{1}%
\end{array}%
\right] &=&-\func{Re}\int \left( i\hat{\psi}_{0}\overline{v_{0}}+iv_{1}%
\overline{\psi _{0}}\right) \;dx \\
&=&-\func{Re}\int \left( i\hat{\psi}_{0}\overline{v_{0}}+\overline{iv_{1}%
\overline{\psi _{0}}}\right) \;dx \\
&=&-\func{Re}\int \left( i\hat{\psi}_{0}\overline{v_{0}}-i\psi _{0}\overline{%
v_{1}}\right) \;dx.
\end{eqnarray*}%
Then (\ref{become}) can be written as follows:
\begin{eqnarray*}
\func{Re}\int \left[ \nabla \psi _{0}\overline{\nabla v_{0}}+W^{\prime
}(x,\psi _{0})\overline{v_{0}}\right] dx &=&\omega \func{Re}\int i\hat{\psi}%
_{0}\overline{v_{0}}\;dx \\
\func{Re}\int \hat{\psi}_{0}\overline{v_{1}}\ dx &=&-\omega \func{Re}\int
i\psi _{0}\overline{v_{1}}\;dx.
\end{eqnarray*}%
Then%
\begin{eqnarray}
-\Delta \psi _{0}+W^{\prime }(x,\psi _{0}) &=&i\omega \hat{\psi}_{0}  \notag
\\
\hat{\psi}_{0} &=&-i\omega \psi _{0}  \label{lui}
\end{eqnarray}

So we get (\ref{ella}). From (\ref{ella}) and (\ref{lui}) we easily verify
that (\ref{evol}) solves (\ref{ham1}), (\ref{ham2}).

$\square $

\subsection{Existence results for NKG}

The following Theorem holds:

\begin{theorem}
\label{klein}Assume that $W$ satisfies NKG-i),...NKG-iiii) and that
\begin{equation}
\alpha <h_{0}  \label{il}
\end{equation}%
where $h_{0}$ is defined by (\ref{h0}) and (\ref{h1}). Then equation (NKG)
admits hylomorphic solitons having the following form%
\begin{equation*}
U(t,x)=(\psi _{0}(x)e^{-i\omega t},-i\omega \psi _{0}(x)e^{-i\omega t}).
\end{equation*}
\end{theorem}

In order to prove the existence of hylomorphic solitons, we will use Th. \ref%
{due}. Clearly the energy $E$ and the hylenic charge $C$ have the form (\ref%
{elle}) and (\ref{cu}) respectively, with%
\begin{equation*}
X=\left\{ u=\left( \psi ,\hat{\psi}\right) \in H^{1}(\mathbb{R}^{N},\mathbb{C%
})\times L^{2}(\mathbb{R}^{N},\mathbb{C})\right\}
\end{equation*}%
\begin{equation}
\left\langle Lu,u\right\rangle =\int \left( \left\vert \hat{\psi}\right\vert
^{2}+\left\vert \nabla \psi \right\vert ^{2}+h^{2}\left\vert \psi
\right\vert ^{2}\right) dx;\ K(u)=\int N(\psi )dx,
\end{equation}

\begin{equation}
\left\langle L_{0}u,u\right\rangle =C(u)=-\func{Re}\int i\hat{\psi}\overline{%
\psi }\;dx;\ K_{0}(u)=0.
\end{equation}

Now let us compute $e_{0}$ and $\Lambda _{\ast }$ defined by (\ref{brutta})
and (\ref{brutta1}).

\begin{lemma}
\label{ve}Assume that $W$ satisfies NKG-i,...NKG-iiii), then
\begin{eqnarray}
e_{0} &\geq &h_{0} \\
\Lambda _{\ast } &\leq &\alpha .
\end{eqnarray}
\end{lemma}

\textbf{Proof}. By (\ref{brutta}) we have
\begin{eqnarray}
e_{0} &=&\ \inf \ \frac{\frac{1}{2}\left\langle Lu,u\right\rangle _{Q}}{%
\left\langle L_{0}u,u\right\rangle _{Q}}=\ \inf \ \frac{\frac{1}{2}%
\int_{Q}\left( \left\vert \hat{\psi}\right\vert ^{2}+\left\vert \nabla \psi
\right\vert ^{2}+h\left( x\right) ^{2}\left\vert \psi \right\vert
^{2}\right) dx}{\left\vert \func{Re}\int_{Q}i\hat{\psi}\overline{\psi }%
\;dx\right\vert }  \label{birba} \\
&\geq &\ \inf \ \ \frac{\frac{1}{2}\int_{Q}\left( \left\vert \hat{\psi}%
\right\vert ^{2}+h_{0}^{2}\left\vert \psi \right\vert ^{2}\right) dx}{%
\int_{Q}\left\vert \hat{\psi}\right\vert \cdot \left\vert \psi \right\vert
\;dx}\geq \ \inf \ \frac{h_{0}\int_{Q}\left\vert \hat{\psi}\right\vert \cdot
\left\vert \psi \right\vert \;dx}{\int_{Q}\left\vert \hat{\psi}\right\vert
\cdot \left\vert \psi \right\vert \;dx}=h_{0}.  \notag
\end{eqnarray}%
Then%
\begin{equation*}
e_{0}\geq h_{0}
\end{equation*}

Let us now prove that%
\begin{equation*}
\Lambda _{\ast }\leq \alpha
\end{equation*}

Let $R>0;$ set
\begin{equation}
u_{R}=\left\{
\begin{array}{cc}
\bar{s} & if\;\;|x|<R \\
0 & if\;\;|x|>R+1 \\
\frac{|x|}{R}\bar{s}-(\left\vert x\right\vert -R)\frac{R+1}{R}\bar{s} &
if\;\;R<|x|<R+1%
\end{array}%
\right.  \label{inff}
\end{equation}
and set $\psi =u_{R},\ $and\ $\hat{\psi}=\alpha u_{R}.$

Then
\begin{align*}
\Lambda _{\ast }& =\underset{\psi ,\hat{\psi}}{\ \inf }\ \frac{\int \left(
\frac{1}{2}\left\vert \hat{\psi}\right\vert ^{2}+\frac{1}{2}\left\vert
\nabla \psi \right\vert ^{2}+W(\psi )\right) dx}{\left\vert \func{Re}\int i%
\hat{\psi}\overline{\psi }\;dx\right\vert } \\
& \leq \frac{\int \left( \frac{1}{2}\alpha ^{2}\left\vert u_{R}\right\vert
^{2}+\frac{1}{2}\left\vert \nabla u_{R}\right\vert ^{2}+W(u_{R})\right) dx}{%
\alpha \int \left\vert u_{R}\right\vert ^{2}\;dx} \\
& \leq \ \frac{\int_{\left\vert x\right\vert <R}\left( \frac{1}{2}\alpha
^{2}\left\vert u_{R}\right\vert ^{2}+W(u_{R})\right) dx}{\alpha
\int_{\left\vert x\right\vert <R}\left\vert u_{R}\right\vert ^{2}\;dx} \\
& +\frac{\int_{R<\left\vert x\right\vert <R+1}\left( \frac{1}{2}\alpha
^{2}\left\vert u_{R}\right\vert ^{2}+\frac{1}{2}\left\vert \nabla
u_{R}\right\vert ^{2}+W(u_{R})\right) dx}{\alpha \int_{\left\vert
x\right\vert <R}\left\vert u_{R}\right\vert ^{2}dx} \\
& =\frac{1}{2}\alpha +\frac{\int_{\left\vert x\right\vert <R}W(\bar{s})dx}{%
\alpha \int_{\left\vert x\right\vert <R}\left\vert \bar{s}\right\vert ^{2}dx}%
+O\left( \frac{1}{R}\right) ;
\end{align*}

Then, by NKG-ii and (\ref{inff})%
\begin{equation*}
\Lambda _{\ast }\leq \frac{1}{2}\alpha +\frac{\int_{\left\vert x\right\vert
<R}\frac{1}{2}\bar{s}^{2}\alpha ^{2}R^{N}}{\alpha \int_{\left\vert
x\right\vert <R}\left\vert \bar{s}\right\vert ^{2}dx}+O\left( \frac{1}{R}%
\right) =\alpha +O\left( \frac{1}{R}\right) .
\end{equation*}%
Then, we get%
\begin{equation*}
\Lambda _{\ast }\leq \alpha
\end{equation*}

$\square $

\textbf{Proof of Theorem} \ref{klein}:

By Lemma \ref{ve} and assumption (\ref{il}) we deduce that $\Lambda _{\ast
}<e_{0}.$ By standard calculations it can be shown that under the
assumptions NKG-i),...,NKG-iiii) the functionals $E$ and $C,$ defined by (%
\ref{en}) and (\ref{cn}), satisfy (E1,2,3,4) of section \ref{gf}. Then, by
using Th. \ref{due} and remark \ref{figo}, we deduce that equation (\ref{KG}%
) admits hylomorphic solitons. Since these solitons are minimizers of the
energy $E$ on the manifold $\left\{ u\in X:C(u)=\sigma \right\} ,$ we easily
get, by Proposition \ref{pipi}, that they are solutions of (\ref{KG}) of the
type $U(t,x)=(\psi _{0}(x)e^{-i\omega t},-i\omega \psi _{0}(x)e^{-i\omega
t}) $ with $\psi _{0},\omega $ satisfying (\ref{ella}).

$\square $

\bigskip \bigskip


\begin{thebibliography}{99}
\bibitem{Ba-Be-R.} \textsc{Badiale M.}, \textsc{Benci V.}, \textsc{Rolando S.%
}, \emph{Solitary waves: physical aspects and mathematical results}, Rend.
Sem. Math. Univ. Pol. Torino \textbf{62} (2004), 107-154.

\bibitem{BBBM} \textsc{Bellazzini J.,Benci V.,Bonanno C., Micheletti A.M.,}
\emph{\ Solitons for the Nonlinear Klein-Gordon-Equation},
(arXiv:0712.1103), to apper on Advanced of Nonlinear Studies (2010).

\bibitem{milano} \textsc{Benci V, }\textit{Hylomorphic solitons, }Milan J.
Math. \textbf{77 }(2009), 271-332.

\bibitem{hylo} \textsc{Bellazzini J., Benci V., Bonanno C., Sinibaldi E.,}
\emph{\ Hylomorphic solitons} \emph{in the nonlinear Klein-Gordon equation,}
Dynamics of Partial Differential Equations, \textbf{6 }(2009), 311-336.

\bibitem{BBGM} \textsc{Bellazzini J.,Benci V., Ghimenti M., Micheletti A.M.,}
\emph{\ On the existence of the fundamental eigenvalue of an elliptic
problem in $\mathbb{R}^{N}$ }, Adv. Nonlinear Stud. \textbf{7} (2007),
439--458

\bibitem{befo} \textsc{Benci V. Fortunato D.,}\textit{\ Solitary waves in
Abelian Gauge Theories, }Adv. Nonlinear Stud. \textbf{3 }(2008), 327-352.

\bibitem{befo08} \textsc{Benci V. Fortunato D.,}\textit{\ Existence of
hylomorphic solitary waves in Klein-Gordon and in Klein-Gordon-Maxwell
equations,} Rend. Lincei Mat. Appl. \textbf{20 }(2009).

\bibitem{befogranas} \textsc{Benci V. Fortunato D.,}\textit{\ \ Solitary
waves in the Nonlinear Wave equation and in Gauge Theories, }Journal of
fixed point theory and Applications, \textbf{1}, n.1 (2007), 61-86.

\bibitem{befo10} \textsc{Benci V. Fortunato D., }\textit{Spinning Q-balls
for the Klein-Gordon-Maxwell Equations, }Commun. Math. Phys., \textbf{295}
(2010), 639-668.

\bibitem{bonanno} \textsc{Bonanno C.,}\textit{\ Existence and multiplicity
of stable bound states for the nonlinear Klein-Gordon equation. }Nonlinear
Analysis \textbf{72} (2010), 20-31.

\bibitem{Beres-Lions} \textsc{Berestycki H.}, \textsc{Lions P.L.}, \emph{%
Nonlinear scalar field equations, I - Existence of a ground state}, Arch.
Rational Mech. Anal. \textbf{82} (1983), 313-345.

\bibitem{CL82} \textsc{T.~Cazenave and P.L. Lions,} \emph{Orbital stability
of standing waves for some nonlinear {S}chr\"{o}dinger equations}, Comm.
Math. Phys. \textbf{85} (1982), no.~4, 549--561.

\bibitem{Coleman86} \textsc{Coleman S.,}\emph{Q-Balls}, Nucl. Phys. B262
(1985) 263-283; erratum: B269 (1986)\ 744-745.

\bibitem{Gelfand} \textsc{Gelfand I.M., Fomin S.V.}, \textit{Calculus of
Variations}, Prentice-Hall, Englewood Cliffs, N.J. 1963.

\bibitem{raj} \textsc{Rajaraman R., }\textit{Solitons and instantons, }%
North-Holland,\textit{\ }Amsterdam 1989.

\bibitem{rosen68} \textsc{Rosen G.}, \emph{Particlelike solutions to
nonlinear complex scalar field theories with positive-definite energy
densities}, J. Math. Phys. \textbf{9} (1968), 996-998.

\bibitem{shatah} \textsc{Shatah J.,} \textit{Stable Standing waves of
Nonlinear Klein-Gordon Equations, }Comm. Math. Phys., 91, (1983), 313-327.

\bibitem{strauss} \textsc{Strauss W.A.}, \textit{Existence of solitary waves
in higher dimensions, }Comm. Math. Phys. \textbf{55 }(1977), 149-162

\bibitem{vil} \textsc{Vilenkin A., Shellard E.P.S., }\textit{Cosmic strings
and other topological defects, }Cambridge monographs on mathematical
physics, 1994.

\bibitem{yangL} \textsc{Yang Y.,} \textit{Solitons in Field Theory and
Nonlinear Analysis, }Springer, New York, Berlin, 2000.\bigskip
\end{thebibliography}
\end{document}